\documentclass[11pt, twoside]{amsart}
\usepackage{amsmath, amsthm, amscd, amsfonts, amssymb, graphicx, color, delarray}
\usepackage[bookmarksnumbered, plainpages, backref]{hyperref}

\begin{document}

\title{The Application Domain of Infinite Matrices on Classical Sequence Spaces}

\author{Murat Kiri\c{s}ci}

\address{[Murat Kiri\c{s}ci] Department of Mathematical Education, Hasan Ali Y\"{u}cel Education Faculty,
Istanbul University, Vefa, 34470, Fatih, Istanbul, Turkey \vskip 0.1cm }
\email{mkirisci@hotmail.com, murat.kirisci@istanbul.edu.tr}

\begin{abstract}
The purpose of this paper is twofold. First, we define the new spaces and investigate some topological and structural properties. Also,
we compute dual spaces of new spaces which are help us in the characterization of
matrix mappings. Second, we give some examples related to new spaces.
\end{abstract}

\maketitle

\pagestyle{plain} \makeatletter
\theoremstyle{plain}
\newtheorem{thm}{Theorem}[section]
\numberwithin{equation}{section}
\numberwithin{figure}{section}  
\theoremstyle{plain}
\newtheorem{pr}[thm]{Proposition}
\theoremstyle{plain}
\newtheorem{exmp}[thm]{Example}
\theoremstyle{plain}
\newtheorem{cor}[thm]{Corollary} 
\theoremstyle{plain}
\newtheorem{defin}[thm]{Definition}
\theoremstyle{plain}
\newtheorem{lem}[thm]{Lemma} 
\theoremstyle{plain}
\newtheorem{rem}[thm]{Remark}
\numberwithin{equation}{section}

\section{Introduction}

It is well known that, the $\omega$ denotes the family of all real (or complex)-valued sequences.
$\omega$ is a linear space and each linear subspace of $\omega$ (with the included addition
and scalar multiplication) is called a \emph{sequence space} such as the spaces $c$, $c_{0}$ and
$\ell_{\infty}$, where $c$, $c_{0}$ and $\ell_{\infty}$ denote the set of all convergent
sequences in fields $\mathbb{R}$ or $\mathbb{C}$, the set of all null sequences and the set
of all bounded sequences, respectively. It is clear that the sets $c$, $c_{0}$ and $\ell_{\infty}$
are the subspaces of the $\omega$. Thus, $c$, $c_{0}$ and $\ell_{\infty}$ equipped with a vector space structure,
from a sequence space. By $bs$ and $cs$, we define the spaces of all bounded and convergent series, respectively.\\

\emph{A coordinate space} (or \emph{$K-$space}) is a vector space of numerical sequences, where addition and scalar multiplication are defined pointwise. That is, a sequence space $X$ with a linear topology is called a $K$-space provided each of the maps $p_{i}:X\rightarrow \mathbb{C}$ defined by $p_{i}(x)=x_{i}$ is continuous for all $i\in \mathbb{N}$. A $K-$space is called an \emph{$FK-$space} provided $X$ is a complete linear metric space. An \emph{$FK-$space} whose topology is normable is called a \emph{$BK-$ space}.\\

Let $X$ be a $BK-$space. Then $X$ is said to have monotone norm if $\|x^{[m]}\|\geq \|x^{[n]}\|$ for $m>n$ and $\|x\|=\sup \|x^{[m]}\|$. The spaces $c_{0}$,  $c$, $\ell_{\infty}$,  $cs$,  $bs$ have monotone norms.\\

If a normed sequence space $X$ contains a sequence $(b_{n})$ with the property that for every $x\in X$ there is unique sequence of scalars $(\alpha_{n})$ such that
\begin{eqnarray*}
\lim_{n\rightarrow\infty}\|x-(\alpha_{0}b_{0}+\alpha_{1}b_{1}+...+\alpha_{n}b_{n})\|=0
\end{eqnarray*}
then $(b_{n})$ is called \emph{Schauder basis} for $X$. The series $\sum\alpha_{k}b_{k}$ which has the sum $x$ is then called the expansion of $x$ with respect to $(b_{n})$, and written as $x=\sum\alpha_{k}b_{k}$. An \emph{$FK-$space} $X$ is said to have $AK$ property, if $\phi \subset X$ and $\{e^{k}\}$ is a basis for $X$, where $e^{k}$ is a sequence whose only non-zero term is a $1$ in $k^{th}$ place for each $k\in \mathbb{N}$ and $\phi=span\{e^{k}\}$, the set of all finitely non-zero sequences.An $FK-$space $X\supset \phi$ is said to have $AB$, if $(x^{[n]})$ is a bounded set in $X$ for each $x\in X$.\\

Let $A=(a_{nk})$ be an infinite matrix of complex numbers $a_{nk}$ and $x=(x_{k})\in \omega$, where $k,n\in\mathbb{N}$.
Then the sequence $Ax$ is called as the $A-$transform of $x$ defined by the usual matrix product.
Hence, we transform the sequence $x$ into the sequence $Ax=\{(Ax)_{n}\}$ where
\begin{eqnarray}\label{equ1}
(Ax)_{n}=\sum_{k}a_{nk}x_{k}
\end{eqnarray}
for each $n\in\mathbb{N}$, provided the series on the right hand side of (\ref{equ1}) converges for each $n\in\mathbb{N}$. Let $X$ and $Y$ be two sequence spaces. If $Ax$ exists and is in $Y$ for every sequence $x=(x_{k})\in X$, then we say that
$A$ defines a matrix mapping from $X$ into $Y$, and we denote it by writing $A :X \rightarrow Y$ if and only if the series on
the right hand side of (\ref{equ1}) converges for each $n\in\mathbb{N}$ and every $x\in X$, and we have $Ax=\{(Ax)_{n}\}_{n\in \mathbb{N}}\in Y$
for all $x\in X$.  A sequence $x$ is said to be $A$-summable to $l$ if $Ax$ converges to $l$ which is called the $A$-limit of $x$. Let $X$ be a sequence space and $A$ be an infinite matrix. The sequence space
\begin{eqnarray}\label{eq0}
X_{A}=\{x=(x_{k})\in\omega:Ax\in X\}
\end{eqnarray}
is called the domain of $A$ in $X$ which is a sequence space.\\

The matrix $\Omega=(a_{nk})$ defined by $a_{nk}=k$, $(1\leq k \leq n)$ and $a_{nk}=0$, $(k>n)$, and the matrix $\Gamma=(b_{nk})$ defined by
by $b_{nk}=1/k$, $(1\leq k \leq n)$ and $b_{nk}=0$, $(k>n)$, respectively, i.e.,

\begin{eqnarray*}
a_{nk}= \left[ \begin{array}{cccccc}
1 & 0 & 0 & 0 & \cdots \\
1 & 2 & 0 & 0 & \cdots  \\
1 & 2 & 3 & 0 & \cdots \\
1 & 2 & 3 & 4 & \cdots \\
\vdots & \vdots &  \vdots&  \vdots& \ddots
\end{array} \right]
\quad \textrm{ and } \quad
b_{nk}= \left[ \begin{array}{cccccc}
1 & 0 & 0 & 0 & \cdots \\
1 & 1/2 & 0 & 0 & \cdots  \\
1 & 1/2 & 1/3 & 0 & \cdots \\
1 & 1/2 & 1/3 & 1/4 & \cdots \\
\vdots & \vdots &  \vdots&  \vdots& \ddots
\end{array} \right]
\end{eqnarray*}

We can give the matrices $\Omega^{-1}=(c_{nk})$ and $\Gamma^{-1}=(d_{nk})$ which are inverse of the above matrices by $c_{nk}=1/n$, $(n=k)$, $c_{nk}=-1/n$, $(n-1=k)$, $c_{nk}=0$, $(other)$ and $d_{nk}=n$, $(n=k)$, $d_{nk}=-n$, $(n-1=k)$, $d_{nk}=0$, $(other)$, respectively, i.e.,

\begin{eqnarray*}
c_{nk}= \left[ \begin{array}{cccccc}
1 & 0 & 0 & 0 & \cdots \\
-1/2 & 1/2 & 0 & 0 & \cdots  \\
0 & -1/3 & 1/3 & 0 & \cdots \\
0 & 0 & -1/4 & 1/4 & \cdots \\
\vdots & \vdots &  \vdots&  \vdots& \ddots
\end{array} \right]
\quad \textrm{ and } \quad
d_{nk}= \left[ \begin{array}{cccccc}
1 & 0 & 0 & 0 & \cdots \\
-2 & 2 & 0 & 0 & \cdots  \\
0 & -3 & 3 & 0 & \cdots \\
0 & 0 & -4 & 4 & \cdots \\
\vdots & \vdots &  \vdots&  \vdots& \ddots
\end{array} \right]
\end{eqnarray*}

Now, we show that the matrices $\Omega$ and $\Gamma$ preserve the limits
on the set of all convergent sequences.\\

\begin{thm}
The matrices $\Omega$ and $\Gamma$ are regular.
\end{thm}

\begin{proof}
Take a sequence $x=(x_{k})$. We must show that
if for $n\rightarrow \infty$ and some $L$, $\lim_{n}|x_{k}-L|\rightarrow 0$,
then, $\lim_{n}|b_{nk}x_{k}-L|\rightarrow 0$, where $b_{nk}$ is $\Gamma$ matrix.
Suppose that for $n\rightarrow \infty$ and some $L$, $\lim_{n}|x_{k}-L|\rightarrow 0$, and choose
$\varepsilon >0$. Then, there exists a positive integer $N$ such that $\lim_{n}|x_{k}-L|<\varepsilon$
for $n\geq N$. Then, for $n\geq N$ and $N \in \mathbb{N}$, $\lim_{n}|b_{nk}x_{k}-L|= \lim_{n}|\sum_{k=1}^{n}(k^{-1}x_{k}-L)|<\varepsilon$. Therefore the matrix $\Gamma$ is regular.\\

Similarly, we can show that the matrix $\Omega$ is regular.
\end{proof}

\section{New Spaces and Topological Properties}

Now, we introduce the new sequence spaces derived by the $\Omega-$ and $\Gamma-$ matrices as follows:
\begin{eqnarray*}
&&\ell_{\infty}({\Omega})=\Big\{x=(x_{k})\in \omega: \Omega x \in \ell_{\infty}\}\\
&&c({\Omega})=\{x=(x_{k})\in \omega: \Omega x \in c\}\\
&&c_{0}({\Omega})=\{x=(x_{k})\in \omega: \Omega x \in c_{0}\}
\end{eqnarray*}
and
\begin{eqnarray*}
&&\ell_{\infty}({\Gamma})=\{x=(x_{k})\in \omega: \Gamma x \in \ell_{\infty}\}\\
&&c({\Gamma})=\{x=(x_{k})\in \omega: \Gamma x \in c\}\\
&&c_{0}({\Gamma})=\{x=(x_{k})\in \omega: \Gamma x \in c_{0}\}
\end{eqnarray*}

Let us define the sequences $u=(u_{n})$ and $v=(v_{n})$, as the $\Omega-$transform and $\Gamma-$transform
 of a sequence $x=(x_{k})$, respectively, that is, for $k,n\in \mathbb{N}$, $u_{n}=\left(\Omega x\right)_{n}=\sum_{k=1}^{n}\left|kx_{k}\right|$
 and $v_{n}=\left(\Gamma x\right)_{n}=\sum_{k=1}^{n}\left|k^{-1}x_{k}\right|$.

\begin{thm}\label{isotheo}
The new bounded, convergent and null sequence spaces are norm isomorphic to the
classical sets consisting of the bounded, convergent and null sequences.
\end{thm}

\begin{proof}
We will show that there is a linear isometry between new bounded, convergent, null sequence spaces and classical bounded,
 convergent and null convergent sequence space.
We consider the transformation defined $\Phi$, from $X(\Omega)$ to $X$ by $x \mapsto u =\Phi x = \sum_{k=1}^{n}\left|kx_{k}\right|$, where $X=\{\ell_{\infty}, c, c_{0}\}$. Then, it is
clear that the equality $\Phi(a+b)=\Phi(a)+\Phi(b)$ is holds. Choose $\lambda \in \mathbb{R}$. Then,
\begin{eqnarray*}
\Phi(\lambda a)=\Phi(\lambda a_{k})&=& \sum_{k=1}^{n}\left|\lambda ka_{k}\right|=\lambda \sum_{k=1}^{n}\left|ka_{k}\right|=\lambda \Phi a.
\end{eqnarray*}
Therefore, we can say that $\Phi$ is linear.\\

Choose a sequence $y=(y_{k})$ in $X(\Omega)$ and define the sequence $x=(x_{k})$ such that $x=(c_{nk}y_{k})$,
where $c_{nk}$ is inverse of $\Omega=(a_{nk})$ matrix. Then,

\begin{eqnarray*}
\|x\|_{\ell_{\infty}({\Omega})}&=& \sup_{k}|a_{nk}x_{k}|=\sup_{k} \left|a_{nk}c_{nk}y_{k}\right|_{\ell_{\infty}}=\|y\|_{\ell_{\infty}}.
\end{eqnarray*}
Therefore, we can say that $\Phi$ is norm preserving.\\

Similarly, we can also show that the another spaces are norm isomorphic to classical sequence spaces.
\end{proof}

\begin{thm}\label{normed}
The new bounded, convergent and null sequence spaces are $BK-$spaces with
the norms defined by $\|x\|_{X({\Omega})}=\|\Omega x\|_{\ell_{\infty}}$ and $\|x\|_{X({\Gamma})}=\|\Gamma x\|_{\ell_{\infty}}$, respectively, where
$X=\{\ell_{\infty}, c, c_{0}\}$.
\end{thm}

\begin{proof}
Take a sequence $x=(x_{k})$ in $X(\Omega)$, where $X=\{\ell_{\infty}, c, c_{0}\}$ and define
$f_{k}(x)=x_{k}$ for all $k\in \mathbb{N}$. Then, we have
\begin{eqnarray*}
\|x\|_{X({\Omega})}=\sup \left\{1|x_{1}|+2|x_{2}|+3|x_{3}|+\cdots +k|x_{k}|+\cdots \right\}
\end{eqnarray*}
Therefore, $k|x_{k}|\leq \|x\|_{X({\Omega})} \Rightarrow |x_{k}|\leq K \|x\|_{X({\Omega})}\Rightarrow |f_{k}(x)|\leq K \|x\|_{X({\Omega})}$.
From this result, we say that $f_{k}$ is a continuous linear functional for each $k$. Then, $X(\Omega)$ is a $BK-$space.\\

In the same idea, we can prove that the space $X(\Gamma)$ is a $BK-$space.
\end{proof}

\begin{rem}
We can give the proof of Theorem \ref{normed} in a different way:
From 4.3.1 of \cite{Wil}, we know that if a sequence space $X$ is $BK-$space with respective norm
and $A$ be a triangular infinite matrix, then the matrix domain $X_{A}$ also $BK-$space with
respective norm.
\end{rem}

\begin{thm}\label{AK}
The spaces $X(\Omega)$ and $X(\Gamma)$ have $AK-$property.
\end{thm}

\begin{thm}\label{mon}
The spaces $X(\Omega)$ and $X(\Gamma)$ have monotone norm.
\end{thm}

Theorem \ref{AK} and \ref{mon} can be proved as Theorem 2.4., Theorem 2.6. of \cite{kirisci}.

\begin{rem}\label{rem2}
Any space with a monotone norm has $AB$(10.3.12 of \cite{Wil}).
\end{rem}

\begin{cor}
The spaces $X(\Omega)$ and $X(\Gamma)$ have $AB$.
\end{cor}

Since the isomorphism $\Phi$ , defined in the proof of Theorem \ref{isotheo} is surjective, the inverse image of the basis of the space $X$ is the basis of the spaces $X(\Omega)$, where $X=\{\ell_{\infty},c,c_{0}\}$. Therefore, we have the following:

\begin{thm} The following statements hold:
\begin{itemize}
\item [(i)] Define a sequence $t^{(k)}:=\{t_{n}^{(k)}\}_{n\in\mathbb{N}}$ of elements of the space $X(\Omega)$ for every fixed $k\in \mathbb{N}$ by
\begin{eqnarray*}
t_{n}^{(k)}= \left\{ \begin{array}{ccl}
(-1)^{n-k}k^{-1}&, & \quad (n-1\leq k \leq n)\\
0&, & \quad (1\leq k \leq n-1) \quad \textit{or} \quad (k>n)
\end{array} \right.
\end{eqnarray*}
Then the sequence $\{t^{(k)}\}_{k\in\mathbb{N}}$ is a basis for the space $X(\Omega)$ and if we choose $E_{k}=(\Omega x)_{k}$ for all $k\in\mathbb{N}$, then any $x\in X(\Omega)$ has a unique representation of the form
\begin{eqnarray*}
x:=\sum_{k}E_{k}t^{(k)}.
\end{eqnarray*}
\item [(ii)] Define a sequence $s^{(k)}:=\{s_{n}^{(k)}\}_{n\in\mathbb{N}}$ of elements of the space $X(\Gamma)$ for every fixed $k\in \mathbb{N}$ by
\begin{eqnarray*}
s_{n}^{(k)}= \left\{ \begin{array}{ccl}
(-1)^{n-k}k&, & \quad (n-1\leq k \leq n)\\
0&, & \quad (1\leq k \leq n-1) \quad\textit{or} \quad (k>n)
\end{array} \right.
\end{eqnarray*}
Then the sequence $\{s^{(k)}\}_{k\in\mathbb{N}}$ is a basis for the space $X(\Gamma)$ and if we choose $F_{k}=(\Gamma x)_{k}$ for all $k\in\mathbb{N}$, then any $x\in X(\Gamma)$ has a unique representation of the form
\begin{eqnarray*}
x:=\sum_{k}F_{k}s^{(k)}.
\end{eqnarray*}
\end{itemize}
\end{thm}

\begin{rem}
If a space has a Schauder basis, then it is separable.
\end{rem}

\begin{cor}
The spaces $X(\Omega)$ and $X(\Gamma)$ are separable.
\end{cor}

In this section, we have defined the new spaces derived by infinite matrices and
examined some structural and topological properties.

\section{Dual Spaces}

In this section, we compute dual spaces of new defined spaces.
The beta-, gamma-duals of new defined spaces will help us in the characterization of the matrix mappings.\\

Let $x$ and $y$ be sequences, $X$ and $Y$ be subsets of $\omega$ and $A=(a_{nk})_{n,k=0}^{\infty}$ be an infinite matrix of complex numbers. We write $xy=(x_{k}y_{k})_{k=0}^{\infty}$, $x^{-1}*Y=\{a\in\omega: ax\in Y\}$ and $M(X,Y)=\bigcap_{x\in X}x^{-1}*Y=\{a\in\omega: ax\in Y ~\textrm{ for all }~ x\in X\}$ for the \emph{multiplier space} of $X$ and $Y$. In the special cases of $Y=\{\ell_{1}, cs, bs\}$, we write $x^{\alpha}=x^{-1}*\ell_{1}$,  $x^{\beta}=x^{-1}*cs$,  $x^{\gamma}=x^{-1}*bs$ and $X^{\alpha}=M(X,\ell_{1})$,  $X^{\beta}=M(X,cs)$,  $X^{\gamma}=M(X,bs)$ for the $\alpha-$dual, $\beta-$dual, $\gamma-$dual of $X$. By $A_{n}=(a_{nk})_{k=0}^{\infty}$ we denote the sequence in the $n-$th row of $A$, and we write $A_{n}(x)=\sum_{k=0}^{\infty}a_{nk}x_{k}$ $n=(0,1,...)$ and $A(x)=(A_{n}(x))_{n=0}^{\infty}$, provided $A_{n}\in x^{\beta}$ for all $n$.\\

Now, we will give some lemmas, which are provides convenience in the
compute of the dual spaces and characterize of matrix transformations.\\

\begin{lem}
Matrix transformations between $BK-$spaces are continuous.
\end{lem}

\begin{lem}\cite[Lemma 5.3]{AB2}\label{mtrxtool0}
Let $X, Y$ be any two sequence spaces. $A\in (X: Y_{T})$ if and only if $TA\in (X:Y)$, where
$A$ an infinite matrix and $T$ a triangle matrix.
\end{lem}

\begin{lem}\cite[Theorem 3.1]{AB}\label{mtrxtool}
We define the $B^{T}=(b_{nk})$ depending on a sequence  $a=(a_{k})\in\omega$ and give the inverse of the triangle matrix $T=(t_{nk})$ by
\begin{eqnarray*}
b_{nk}=\sum_{j=k}^na_{j}v_{jk}
\end{eqnarray*}
for all $k,n\in\mathbb{N}$. Then,
\begin{eqnarray*}
X_{T}^{\beta}=\{a=(a_{k})\in\omega: B^{T}\in(X:c)\}
\end{eqnarray*}
and
\begin{eqnarray*}
X_{T}^{\gamma}=\{a=(a_{k})\in\omega: B^{T}\in(X:\ell_{\infty})\}.
\end{eqnarray*}
\end{lem}

Now, we list the following useful conditions.\\
\begin{eqnarray}\label{eq21}
&&\sup_{n}\sum_{k}|a_{nk}|<\infty \\ \label{eq22}
&&\lim_{n \to \infty}a_{nk}-\alpha_{k}=0\\ \label{eq23}
&&\lim_{n \to \infty}\sum_{k}a_{nk} \quad \textit{exists}\\ \label{eq24}
&&\lim_{n \to \infty}\sum_{k}\left|a_{nk}\right|=\sum_{k}\left|\lim_{n \to \infty}a_{nk}\right| \\ \label{eq27}
&&\lim_{n}a_{nk}=0  \quad \textrm{ for all k }\\ \label{eq11}
&&\sup_{m}\sum_{k}\left|\sum_{n=0}^{m}\right|<\infty \\ \label{eq12}
&&\sum_{n}a_{nk} \quad \textit{convergent for all k}\\ \label{eq13}
&&\sum_{n}\sum_{k}a_{nk} \quad \textit{convergent} \\ \label{eq35}
&&\lim_{n}a_{nk} \quad \textrm{ exists for all k } \\ \label{eq14}
&&\lim_{m}\sum_{k}\left|\sum_{n=m}^{\infty}a_{nk}\right|=0
\end{eqnarray}

\begin{lem}\label{lemMTRX}
For the characterization of the class $(X:Y)$ with
$X=\{c_{0}, c, \ell_{\infty}\}$ and $Y=\{\ell_{\infty}, c, cs, bs\}$, we can give the necessary and sufficient
conditions from Table 1, where

\begin{center}
\begin{tabular}{|l | l | l | l |}
\hline \textbf{1.} (\ref{eq21}) & \textbf{2.} (\ref{eq21}), (\ref{eq35}) & \textbf{3.}  (\ref{eq11}) & \textbf{4.} (\ref{eq11}), (\ref{eq12})  \\
\hline \textbf{5.} (\ref{eq21}), (\ref{eq35}), (\ref{eq23}) &  \textbf{6.} (\ref{eq11}), (\ref{eq12}), (\ref{eq13}) & \textbf{7.} (\ref{eq35}), (\ref{eq24}) & \textbf{8.} (\ref{eq14})\\
\hline
\end{tabular}
\end{center}
\end{lem}

\begin{center}
\begin{tabular}{|c | c c c c|}
\hline
To $\rightarrow$ & $\ell_{\infty}$ & $c$ & bs & cs \\ \hline
From $\downarrow$ &  &  & &\\ \hline
$c_{0}$ & \textbf{1.} & \textbf{2.} & \textbf{3.} & \textbf{4.}\\
$c$ & \textbf{1.} & \textbf{5.} & \textbf{3.} & \textbf{6.}\\
$\ell_{\infty}$ & \textbf{1.} & \textbf{7.} & \textbf{3.} & \textbf{8.}\\
\hline
\end{tabular}

\vspace{0.1cm}Table 1\\

\end{center}

For using in the proof of Theorem \ref{dualthm}, we define the matrices $U=(u_{nk})$ and $V=(v_{nk})$ as below:
\begin{eqnarray}\label{mtrxdual}
u_{nk}= \left\{ \begin{array}{ccl}
\frac{a_{k}}{k}-\frac{a_{k+1}}{k+1}&, & \quad (k < n)\\
\frac{a_{n}}{n}&, & \quad (k = n)\\
0&, & \quad (k>n)
\end{array} \right.
\end{eqnarray}

\begin{eqnarray}\label{mtrxdual0}
v_{nk}== \left\{ \begin{array}{ccl}
ka_{k}-(k+1)a_{k+1}&, & \quad (k < n)\\
na_{n}&, & \quad (k = n)\\
0&, & \quad (k>n)
\end{array} \right.
\end{eqnarray}

\begin{thm}\label{dualthm}
The $\beta-$ and $\gamma-$ duals of the new sequence spaces defined by
\begin{eqnarray*}
&&\left[c_{0}(\Omega)\right]^{\beta}= \left\{ a=(a_{k})\in \omega: U\in (c_{0}:c) \right\}\\
&&\left[c(\Omega)\right]^{\beta}= \left\{a=(a_{k})\in \omega: U\in (c:c) \right\}\\
&&\left[\ell_{\infty}(\Omega)\right]^{\beta}= \left\{ a=(a_{k})\in \omega: U\in (\ell_{\infty}:c) \right\}\\
&&\left[c_{0}(\Omega)\right]^{\gamma}= \left\{a=(a_{k})\in \omega: U\in (c_{0}:\ell_{\infty}) \right\}\\
&&\left[c(\Omega)\right]^{\gamma}= \left\{ a=(a_{k})\in \omega: U\in (c:\ell_{\infty}) \right\}\\
&&\left[\ell_{\infty}(\Omega)\right]^{\gamma}= \left\{a=(a_{k})\in \omega: U\in (\ell_{\infty}:\ell_{\infty}) \right\}
\end{eqnarray*}
\end{thm}

\begin{proof}
We will only show the $\beta-$ and $\gamma-$ duals of the new null convergent sequence spaces.
Let $a=(a_{k})\in \omega$. We begin the equality
\begin{eqnarray}\label{dual1}
\sum_{k=1}^{n}a_{k}x_{k}&=&\sum_{k=1}^{n}a_{k}k^{-1}(y_{k}-y_{k-1})
\end{eqnarray}
\begin{eqnarray*}
&=&\sum_{k=1}^{n-1}\left(\frac{a_{k}}{k}-\frac{a_{k+1}}{k+1}\right)y_{k}+\frac{a_{n}}{n}y_{n}=\left(Uy\right)_{n}
\end{eqnarray*}
where $U=(u_{nk})$ is defined by (\ref{mtrxdual}). Using (\ref{dual1}), we can see that
$ax=(a_{k}x_{k})\in cs$ or $bs$ whenever $x=(x_{k})\in c_{0}(\Omega)$ if and only if $Uy\in c$ or $\ell_{\infty}$
whenever $y=(y_{k})\in c_{0}$. Then, from Lemma \ref{mtrxtool0} and Lemma \ref{mtrxtool}, we obtain  the result that
$a=(a_{k})\in \left(c_{0}(\Omega)\right)^{\beta}$ or $a=(a_{k})\in \left(c_{0}(\Omega)\right)^{\gamma}$ if and only if
$U\in (c_{0}:c)$ or $U\in (c_{0}:\ell_{\infty})$, which is what we wished to prove.
\end{proof}

\begin{thm}\label{dualgamma}
The $\beta-$ and $\gamma-$ duals of the new sequence spaces defined by
\begin{eqnarray*}
&&\left[c_{0}(\Gamma)\right]^{\beta}= \left\{ a=(a_{k})\in \omega: V\in (c_{0}:c) \right\}\\
&&\left[c(\Gamma)\right]^{\beta}= \left\{a=(a_{k})\in \omega: V\in (c:c) \right\}\\
&&\left[\ell_{\infty}(\Gamma)\right]^{\beta}= \left\{ a=(a_{k})\in \omega: V\in (\ell_{\infty}:c) \right\}\\
&&\left[c_{0}(\Gamma)\right]^{\gamma}= \left\{a=(a_{k})\in \omega: V\in (c_{0}:\ell_{\infty}) \right\}\\
&&\left[c(\Gamma)\right]^{\gamma}= \left\{ a=(a_{k})\in \omega: V\in (c:\ell_{\infty}) \right\}\\
&&\left[\ell_{\infty}(\Gamma)\right]^{\gamma}= \left\{a=(a_{k})\in \omega: V\in (\ell_{\infty}:\ell_{\infty}) \right\}
\end{eqnarray*}
where $V=(v_{nk})$ is defined by (\ref{mtrxdual0}).
\end{thm}

\section{Matrix Mapping}

Let $X$ and $Y$ be arbitrary subsets of $\omega$. We shall show that,
the characterizations of the classes $(X, Y_{T})$ and $(X_{T},Y)$ can be reduced to that of
$(X, Y)$, where $T$ is a triangle.\\

It is well known that if $h_{c_{0}}(\Delta^{(m)}) \cong c_{0}$, then the equivalence
\begin{eqnarray*}
x\in h_{c_{0}}(\Delta^{(m)}) \Leftrightarrow y\in c_{0}
\end{eqnarray*}
holds. Then, the following theorems will be proved and given
some corollaries which can be obtained to that of Theorems
\ref{mtrxtr1} and \ref{mtrxtr2}. Then, using Lemma \ref{mtrxtool0} and Lemma \ref{mtrxtool}, we have:

\begin{thm}\label{mtrxtr1}
Consider the infinite matrices $A=(a_{nk})$ and $D=(d_{nk})$. These matrices
get associated with each other the following relations:
\begin{eqnarray}\label{eq1}
d_{nk}=\frac{a_{nk}}{k}-\frac{a_{n,k+1}}{k+1}
\end{eqnarray}
for all $k,m, n\in \mathbb{N}$. Then, the following statements are true:\\
\textbf{i.} $A \in (c_{0}(\Omega):Y)$ if and only if $\{a_{nk}\}_{k\in\mathbb{N}} \in [c_{0}(\Omega)]^{\beta}$
for all $n\in \mathbb{N}$ and $D\in (c_{0}:Y)$, where $Y$ be any sequence space.\\
\textbf{ii.}  $A \in (c(\Omega):Y)$ if and only if $\{a_{nk}\}_{k\in\mathbb{N}} \in [c(\Omega)]^{\beta}$
for all $n\in \mathbb{N}$ and $D\in (c:Y)$, where $Y$ be any sequence space.\\
\textbf{iii.}  $A \in (\ell_{\infty}(\Omega):Y)$ if and only if $\{a_{nk}\}_{k\in\mathbb{N}} \in [\ell_{\infty}(\Omega)]^{\beta}$
for all $n\in \mathbb{N}$ and $D\in (\ell_{\infty}:Y)$, where $Y$ be any sequence space.\\
\end{thm}

\begin{proof}
We assume that the (\ref{eq1}) holds between the entries of $A=(a_{nk})$ and $D=(d_{nk})$.
Let us remember that from Theorem \ref{isotheo}, the spaces $c_{0}(\Omega)$ and $c_{0}$ are linearly isomorphic. Firstly,
we choose any $y=(y_{k})\in c_{0}$ and consider $A \in (c_{0}(\Omega):Y)$. Then, we obtain that $D\Omega$ exists and
$\{a_{nk}\}\in (c_{0}(\Omega))^{\beta}$ for all $k\in \mathbb{N}$. Therefore, the necessity of (\ref{eq1}) yields and
$\{d_{nk}\}\in c_{0}^{\beta}$ for all $k,n\in \mathbb{N}$. Hence, $Dy$ exists for each $y\in c_{0}$. Thus, if we take
$m\rightarrow \infty$ in the equality
\begin{eqnarray*}
\sum_{k=1}^{m}a_{nk}x_{k}=\sum_{k=1}^{m}a_{nk}\left(\frac{a_{nk}}{k}-\frac{a_{n,k+1}}{k+1}\right)y_{k}
\end{eqnarray*}
for all $m,n\in \mathbb{N}$, then, we understand that $Dy=Ax$. So, we obtain that $D\in (c_{0}:Y)$.\\

Now, we consider that $\{a_{nk}\}_{k\in\mathbb{N}} \in (c_{0}(\Omega))^{\beta}$
for all $n\in \mathbb{N}$ and $D\in (c_{0}:Y)$. We take any $x=(x_{k})\in c_{0}(\Omega)$. Then, we can see that
$Ax$ exists. Therefore, from the equality

\begin{eqnarray*}
\sum_{k}d_{nk}y_{k}=\sum_{k}a_{nk}x_{k}
\end{eqnarray*}
for all $n\in \mathbb{N}$, we obtain that $Ax=Dy$. Therefore, this shows that $A \in (c_{0}(\Omega):Y)$.
\end{proof}

\begin{thm}\label{mtrxtr2}
Consider that the infinite matrices $A=(a_{nk})$
and $E=(e_{nk})$ with
\begin{eqnarray}\label{eq2}
e_{nk}:=\sum_{k=1}^{\infty}\sum_{j=1}^{n}a_{jk}.
\end{eqnarray}
Then, the following statements are true:\\
\textbf{i.} $A=(a_{nk})\in (X:c_{0}(\Omega))$ if and only if $E\in (X :c_{0})$\\
\textbf{ii.} $A=(a_{nk})\in (X:c(\Omega))$ if and only if $E\in (X :c)$\\
\textbf{iii.} $A=(a_{nk})\in (X:\ell_{\infty}(\Omega))$ if and only if $E\in (X :\ell_{\infty})$
\end{thm}

\begin{proof}

We take any $z=(z_{k})\in X$. Using the (\ref{eq2}), we have

\begin{eqnarray}\label{trans1}
\sum_{k=1}^{m}e_{nk}z_{k}=\sum_{k=1}^{m}\left(\sum_{k=1}^{\infty}\sum_{j=1}^{m}jb_{jk}\right)z_{k}
\end{eqnarray}

for all $m,n\in \mathbb{N}$. Then, for $m\rightarrow \infty$, equation (\ref{trans1}) gives us that $(Ez)_{n}=\{\Omega(Az)\}_{n}$.
Therefore, one can immediately observe from this that $Az\in c_{0}(\Omega)$ whenever $z\in X$ if and only if $Ez\in c_{0}$ whenever $z\in X$. Thus, the proof is completed.
\end{proof}

\begin{thm}
Consider the infinite matrices $A=(a_{nk})$ and $F=(f_{nk})$. These matrices
get associated with each other the following relations:
\begin{eqnarray}
f_{nk}=ka_{nk}-(k+1)a_{n,k+1}
\end{eqnarray}
for all $k,m, n\in \mathbb{N}$. Then, the following statements are true:\\
\textbf{i.} $A \in (c_{0}(\Gamma):Y)$ if and only if $\{a_{nk}\}_{k\in\mathbb{N}} \in [c_{0}(\Gamma)]^{\beta}$
for all $n\in \mathbb{N}$ and $F\in (c_{0}:Y)$, where $Y$ be any sequence space.\\
\textbf{ii.}  $A \in (c(\Gamma):Y)$ if and only if $\{a_{nk}\}_{k\in\mathbb{N}} \in [c(\Gamma)]^{\beta}$
for all $n\in \mathbb{N}$ and $F\in (c:Y)$, where $Y$ be any sequence space.\\
\textbf{iii.}  $A \in (\ell_{\infty}(\Gamma):Y)$ if and only if $\{a_{nk}\}_{k\in\mathbb{N}} \in [\ell_{\infty}(\Gamma)]^{\beta}$
for all $n\in \mathbb{N}$ and $F\in (\ell_{\infty}:Y)$, where $Y$ be any sequence space.\\
\end{thm}

\begin{thm}
Consider that the infinite matrices $A=(a_{nk})$
and $G=(g_{nk})$ with
\begin{eqnarray}
g_{nk}:=\sum_{k=1}^{\infty}\sum_{j=1}^{n}j^{-1}a_{jk}
\end{eqnarray}
Then, the following statements are true:\\
\textbf{i.} $A=(a_{nk})\in (X:c_{0}(\Gamma))$ if and only if $G\in (X :c_{0})$\\
\textbf{ii.} $A=(a_{nk})\in (X:c(\Gamma))$ if and only if $G\in (X :c)$\\
\textbf{iii.} $A=(a_{nk})\in (X:\ell_{\infty}(\Gamma))$ if and only if $G\in (X :\ell_{\infty})$
\end{thm}

\section{Examples}

If we choose any sequence spaces $X$ and $Y$ in Theorem \ref{mtrxtr1} and \ref{mtrxtr2} in previous section,
then, we can find several consequences in every choice. For example, if we take
the space $\ell_{\infty}$ and the spaces which are isomorphic to $\ell_{\infty}$
instead of $Y$ in Theorem \ref{mtrxtr1}, we obtain the following examples:

\begin{exmp}\label{exmpE}
The Euler sequence space $e_{\infty}^{r}$ is defined by $e_{\infty}^{r}=\{x\in \omega: \sup_{n\in\mathbb{N}}|\sum_{k=0}^{n}\binom{n}{k}(1-r)^{n-k}r^{k}x_{k}|<\infty\}$ (\cite{BF2} and \cite{BFM}).
We consider the infinite matrix $A=(a_{nk})$ and define the matrix $H=(h_{nk})$ by
\begin{eqnarray*}
h_{nk}=\sum_{j=0}^{n}\binom{n}{j}(1-r)^{n-j}r^{j}a_{jk}  \quad \quad (k,n\in \mathbb{N}).
\end{eqnarray*}
If we want to get necessary and sufficient conditions for the class $(c_{0}(\Omega): e_{\infty}^{r})$ in Theorem \ref{mtrxtr1},
then, we replace the entries of the matrix $A$ by those of the matrix $H$.
\end{exmp}

\begin{exmp}\label{exmpR}
Let $T_{n}=\sum_{k=0}^{n}t_{k}$ and $A=(a_{nk})$ be an infinite matrix. We define the matrix $P=(p_{nk})$ by
\begin{eqnarray*}
p_{nk}=\frac{1}{T_{n}}\sum_{j=0}^{n}t_{j}a_{jk}  \quad \quad (k,n\in \mathbb{N}).
\end{eqnarray*}
Then, the necessary and sufficient conditions in order for $A$ belongs to the class $(c_{0}(\Omega):r_{\infty}^{t})$
are obtained from in Theorem \ref{mtrxtr1} by replacing the entries of the matrix $A$ by those of the matrix $P$;
 where $r_{\infty}^{t}$ is the space of all sequences whose $R^{t}-$transforms is in the space $\ell_{\infty}$ \cite{malk}.
\end{exmp}

\begin{exmp}
In the space $r_{\infty}^{t}$, if we take $t=e$, then, this space become to the Cesaro sequence space of non-absolute type $X_{\infty}$ \cite{NgLee}.
As a special case, Example \ref{exmpR} includes the characterization of class $((c_{0}(\Omega):r_{\infty}^{t})$.
\end{exmp}

\begin{exmp}
The Taylor sequence space $t_{\infty}^{r}$ is defined by $t_{\infty}^{r}=\{x\in \omega: \sup_{n\in\mathbb{N}}|\sum_{k=n}^{\infty}\binom{k}{n}(1-r)^{n+1}r^{k-n}x_{k}|<\infty\}$ (\cite{kirisci4}).
We consider the infinite matrix $A=(a_{nk})$ and define the matrix $T=(t_{nk})$ by
\begin{eqnarray*}
t_{nk}=\sum_{k=n}^{\infty}\binom{k}{n}(1-r)^{n+1}r^{k-n}a_{jk}  \quad \quad (k,n\in \mathbb{N}).
\end{eqnarray*}
If we want to get necessary and sufficient conditions for the class $(c_{0}(\Omega): t_{\infty}^{r})$ in Theorem \ref{mtrxtr1},
then, we replace the entries of the matrix $A$ by those of the matrix $T$.
\end{exmp}

If we take the spaces $c$, $cs$ and $bs$ instead of $X$ in Theorem \ref{mtrxtr2}, or $Y$ in Theorem \ref{mtrxtr1}
we can write the following examples. Firstly, we give some conditions and following lemmas:

\begin{eqnarray}\label{eq27x}
&&\lim_{k}a_{nk}=0  \quad \textrm{ for all n }, \\ \label{eq28}
&&\lim_{n \to \infty}\sum_{k}a_{nk}=0, \\ \label{eq26x}
&&\lim_{n \to \infty}\sum_{k} |a_{nk}|=0, \\ \label{eq33}
&&\lim_{n \to \infty}\sum_{k} |a_{nk}-a_{n,k+1}|=0, \\ \label{eq34}
&&\sup_{n}\sum_{k}\left|a_{nk}-a_{n,k+1}\right|<\infty \\ \label{eq36}
&&\lim_{k}\left(a_{nk}-a_{n,k+1}\right) \textrm{ exists for all k } \\ \label{eq37}
&&\lim_{n \to \infty}\sum_{k}\left|a_{nk}-a_{n,k+1}\right|=\sum_{k}\left|\lim_{n \to \infty}(a_{nk}-a_{n,k+1})\right|\\ \label{eq38}
&&\sup_{n}\left|\lim_{k}a_{nk}\right|<\infty
\end{eqnarray}

\begin{lem}
Consider that the $X\in\{\ell_{\infty}, c, bs, cs\}$ and $Y\in \{c_{0}\}$.
The necessary and sufficient conditions for $A\in (X:Y)$ can be read the following, from Table 2:

\begin{center}
\begin{tabular}{|l | l | l | l |}
\hline \textbf{9.} (\ref{eq26x}) & \textbf{10.} (\ref{eq21}), (\ref{eq27}), (\ref{eq28}) & \textbf{11.} (\ref{eq27x}), (\ref{eq33}) & \textbf{12.} (\ref{eq27}), (\ref{eq34})  \\
\hline \textbf{13.} (\ref{eq27x}), (\ref{eq36}), (\ref{eq37}) &  \textbf{14.} (\ref{eq34}), (\ref{eq35})  & \textbf{15.} (\ref{eq27x}), (\ref{eq34}) & \textbf{16.} (\ref{eq34}), (\ref{eq38}) \\ \hline
\end{tabular}
\end{center}

\end{lem}

\begin{center}
\begin{tabular}{|c | c c  c c|}
\hline
From $\rightarrow$ & $\ell_{\infty}$ & $c$ &  $bs$ & $cs$ \\ \hline
To $\downarrow$ &    & & &\\ \hline
$c_{0}$ & \textbf{9.} & \textbf{10.} & \textbf{11.} & \textbf{12.}\\
$c$ & \textbf{7.} & \textbf{5.} & \textbf{13.} & \textbf{14.}\\
$\ell_{\infty}$ & \textbf{1.} & \textbf{1.} & \textbf{15.} & \textbf{16.}\\
\hline
\end{tabular}

\vspace{0.1cm}Table 2\\
\end{center}

\begin{exmp}\label{ExpOmega1}
We choose $X\in \{c_{0}(\Omega), c(\Omega), \ell_{\infty}(\Omega)\}$ and $Y\in \{\ell_{\infty}, c, cs, bs, f\}$.
The necessary and sufficient conditions for $A\in (X:Y)$ can be taken from the Table 3:
\end{exmp}
\begin{itemize}
  \item[\textbf{1a.}] (\ref{eq21}) holds with $d_{nk}$ instead of $a_{nk}$.
 \item[\textbf{2a.}] (\ref{eq21}), (\ref{eq35}) hold with $d_{nk}$ instead of $a_{nk}$.
\item[\textbf{3a.}] (\ref{eq11}) holds with $d_{nk}$ instead of $a_{nk}$.
\item[\textbf{4a.}] (\ref{eq11}), (\ref{eq12}) hold with $d_{nk}$ instead of $a_{nk}$.
\item[\textbf{5a.}] (\ref{eq21}), (\ref{eq35}), (\ref{eq23}) hold with $d_{nk}$ instead of $a_{nk}$.
\item[\textbf{6a.}] (\ref{eq11}), (\ref{eq12}), (\ref{eq13}) hold with $d_{nk}$ instead of $a_{nk}$.
\item[\textbf{7a.}] (\ref{eq35}), (\ref{eq24}) hold with $d_{nk}$ instead of $a_{nk}$.
\item[\textbf{8a.}] (\ref{eq14}) holds with $d_{nk}$ instead of $a_{nk}$.
\end{itemize}

\begin{center}
\begin{tabular}{|c | c c c c |}
\hline
To $\rightarrow$ & $\ell_{\infty}$ & $c$ & $bs$ & $cs$ \\ \hline
From $\downarrow$ &  &  & & \\ \hline
$c_{0}(\Omega)$ & \textbf{1a.} & \textbf{2a.} & \textbf{3a.} & \textbf{4a.} \\
$c(\Omega)$ & \textbf{1a.} & \textbf{5a.} & \textbf{3a.} & \textbf{6a.} \\
$\ell_{\infty}(\Omega)$ & \textbf{1a.} & \textbf{7a.} & \textbf{3a.} & \textbf{8a.}\\
\hline
\end{tabular}

\vspace{0.1cm}Table 3\\
\end{center}

\begin{exmp}\label{ExpOmega2}
Consider that the $X\in\{\ell_{\infty}, c, bs, cs\}$ and $Y\in \{c_{0}(\Omega), c(\Omega), \ell_{\infty}(\Omega)\}$.
The necessary and sufficient conditions for $A\in (X:Y)$ can be read the following, from Table 4:
\begin{itemize}
\item[\textbf{9a.}] (\ref{eq26x}) holds with $e_{nk}$ instead of $a_{nk}$.
\item[\textbf{10a.}] (\ref{eq21}), (\ref{eq27}), (\ref{eq28}) hold with $e_{nk}$ instead of $a_{nk}$.
\item[\textbf{11a.}] (\ref{eq27x}), (\ref{eq33}) hold with $e_{nk}$ instead of $a_{nk}$.
\item[\textbf{12a.}] (\ref{eq27}), (\ref{eq34}) hold with $e_{nk}$ instead of $a_{nk}$.
\item[\textbf{13a.}] (\ref{eq27x}), (\ref{eq36}), (\ref{eq37}) hold with $e_{nk}$ instead of $a_{nk}$.
\item[\textbf{14a.}] (\ref{eq34}), (\ref{eq35}) hold with $e_{nk}$ instead of $a_{nk}$.
\item[\textbf{15a.}] (\ref{eq27x}), (\ref{eq34}) hold with $e_{nk}$ instead of $a_{nk}$.
\item[\textbf{16a.}] (\ref{eq34}), (\ref{eq38}) hold with $e_{nk}$ instead of $a_{nk}$.
\end{itemize}
\end{exmp}

\begin{center}
\begin{tabular}{|c | c c  c c|}
\hline
From $\rightarrow$ & $\ell_{\infty}$ & $c$ &  $bs$ & $cs$ \\ \hline
To $\downarrow$ &    & & &\\ \hline
$c_{0}(\Omega)$ & \textbf{9a.} & \textbf{10a.} & \textbf{11a.} & \textbf{12a.}\\
$c(\Omega)$ & \textbf{7a.} & \textbf{5a.} & \textbf{13a.} & \textbf{14a.} \\
$\ell_{\infty}(\Omega)$ & \textbf{1a.} & \textbf{1a.} & \textbf{15a.} & \textbf{16a.}\\
\hline
\end{tabular}

\vspace{0.1cm}Table 4\\
\end{center}

With the same idea of Example \ref{ExpOmega1} and Example \ref{ExpOmega2}, we can write the examples related to the $\Gamma$ matrix as table form.
In examples which are writing with $\Gamma$ matrix, we use the $f_{nk}$ and $g_{nk}$.

\section{Conclusion}
We know that the most general linear operators between two sequence spaces is given by an infinite matrix.
he theory of matrix transformations deals with establishing necessary and sufficient conditions on the
   entries of a matrix to map a sequence space $X$ into a sequence space $Y$. This is a natural generalization
   of the problem to characterize all summability methods given by infinite matrices that preserve convergence.\\

In this work, we construct a new sequence spaces by means of the matrix domain with two infinite matrices.
We examine to some properties such as isomorphism, $BK-$space, $AK-$ and $AB-$properties, monotone norm.
Also, we give dual spaces and later the necessary and sufficient conditions on
the matrix transformations of the classes $(X:Y)$. Afterward, in the last section, we obtain several examples related to new spaces.

\end{document}